\documentclass[11pt]{amsart}
 \usepackage{amsmath,amssymb,enumerate}
 \usepackage[all]{xy}
 \usepackage{xcolor}

 \newtheorem{pro}{Proposition}[section]
 \newtheorem{cor}[pro]{Corollary}
 \newtheorem{lem}[pro]{Lemma}

  \newtheorem{defi}[pro]{Definition}
 
 \newtheorem{thm}[pro]{Theorem}

 \newcommand{\C}{\mathbb C}

   \newcommand{\G}{\mathbb G}
\usepackage[colorlinks=true,
            linkcolor=blue,
            urlcolor=blue,
            citecolor=blue]{hyperref}

 \numberwithin{equation}{section}
\usepackage{hyperref}
\begin{document}
\title[Logarithmic Potentials]{Logarithmic Potentials on $\mathbb P^n$}
\setcounter{tocdepth}{1}
\author{  Fatima Zahra Assila}

\date{\today} 

 \address{Ibn Tofail university\\ faculty of sciences\\ PO 242 Kenitra Morocco}
\email{fatima.zahra.assila@uit.ac.ma}

\maketitle

\begin{abstract}
%\selectlanguage{english}
We study the projective logarithmic potential   $\mathbb{G}_{\mu}$ of a Probability  measure $\mu$ on the complex projective space $\mathbb{P}^{n}$. We prove that the Range of the  operator $\mu\longrightarrow \mathbb{G}_{\mu}$  is contained in  the (local) domain of definition of the complex Monge-Amp\`ere operator acting on the class of quasi-plurisubharmonic functions on $\mathbb{P}^n$ with respect to the Fubini-Study metric. Moreover, when the measure $\mu $ has no atom,
 we show that the complex Monge-Amp\`ere measure of its Logarithmic potential is an absolutely continuous measure with respect to the Fubini-Study volume form on $\mathbb{P}^{n}$ 
\end{abstract}

\section{Introduction and statement of the results}
Logarithmic potentials of Borel measures in the complex plane play a  fundamental role in Logarithmic Potential Theory.  This du to the fact that this theory is associated to the Laplace operator which is a linear elliptic partial differential operator of second order. It is well known that in higher dimension plurisubharmonic functions are rather connected to the complex Monge-Amp\`ere operator which is a fully non-linear second order partial differential operator. Therefore Pluripotential theory cannot be described by logarithmic potential. However the class of logarithmic potentials gives a nice class of plurisubharmonic functions which turns out to be in the local domain of definition of the complex Monge-Amp\`ere operator. This study was carried out by Carlehed \cite{5} in the case of a compactly supported measures on $\mathbb{C}^{n}$ or a bounded hyperconvex domain in $\mathbb{C}^{n}$.

 Our main goal is to extend this study to the complex projective space motivated by the fact that the complex Monge-Amp\`ere operator plays an important role in K\"ahler geometry (see \cite{13}). A large class of singular potentials on which the complex Monge-Amp\`ere is well defined was introduced (see  \cite{12}, \cite{8}, \cite{4}). However the global domain of definition of the complex Monge-Amp\`ere operator on compact K\"ahler manifolds is not yet well understood. Using  the characterization of the local domain of definition  given by Cegrell and Blocki (see \cite{2}, \cite{3}, \cite{7}), we show  that it is contained in the local domain of definition of the complex Monge-Amp\`ere operator on the complex projective space $(\mathbb{P}^{n},\omega)$ equipped with the Fubini-Study metric. 
 
 Let $\mu$ be a probability measure on $\mathbb{P}^{n}$. 
Then its projective logarithmic potential is defined on $\mathbb{P}^{n}$ as follows :
$$
\mathbb{G}_{\mu}(\zeta):=\int_{\mathbb{P}^{n}} G(\zeta,\eta)d\mu(\eta)\quad\hbox{where}\quad G(\zeta,\eta):=\log{\vert\zeta\wedge\eta\vert\over\vert\zeta\vert\vert\eta\vert}
$$
\smallskip

\begin{thm}\label{thm1}
Let $\mu$ be a probability measure on $\mathbb{P}^{n}$. Then the following properties hold.\\
1. The potential $\mathbb{G}_{\mu}$ is a negative $\omega$-plurisubharmonic function on $\mathbb{P}^{n}$ normalized by the following condition : 
$$
\int_{\mathbb{P}^{n}} \mathbb{G}_{\mu} \, \omega_{FS}^n = -\alpha_{n},
$$
where $\alpha_{n}$ is a numerical constant.\\
2. $\mathbb{G}_{\mu}\in W^{1,p}(\mathbb{P}^{n})$ for any $0 < p < 2n$.\\
3. $\mathbb{G}_{\mu}\in DMA_{loc}(\mathbb{P}^{n},\omega)$.
\end{thm}

\smallskip

We also show a regularizing property of the operator $\mu\rightarrow\mathbb{G}_{\mu}$ acting on probability measures on $\mathbb{P}^{n}$.
\begin{thm}\label{thm2}
Let $\mu$ be a probability measure on $\mathbb{P}^{n}$ with no atoms.  Then the Monge-Amp\`ere measure $(\omega +dd^{c}\mathbb{G}_{\mu})^{n}$  is absolutely continuous with respect to the Fubini-Study volume form on $\mathbb{P}^{n}$.
\end{thm}

\section{The logarithmic potential, proof of Theorem 1.1}

The complex projective space can be covered by a finite number of charts given by $ \mathcal{U}_{k}:=\{[\zeta_{0},\zeta_{1},\cdots,\zeta_{n}]\in\mathbb{P}^{n}\ :\ \zeta_{k}\not=0\}\ ( 0\leq k\leq n)$ and the corresponding coordinate chart is defined on $\mathcal{U}_{k}$ by the formula

$$
z^{k}(\zeta)=z^{k}:=(z_{j}^{k})_{0\leq j\leq n, j\not=k}\quad\hbox{where}\quad z_{j}^{k}:={\zeta_{j}\over\zeta_{k}}\quad\hbox{for}\quad j\not=k
$$
 The Fubini-Study metric $\omega = \omega_{FS}$ is given on $\mathcal U_k$  by  $\omega|_{\mathcal{U}_{k}}={1\over 2}dd^{c}\log(1+\vert z^{k}\vert^{2})$. The projective logaritmic kernel on $\mathbb{P}^{n}\times\mathbb{P}^{n}$ is naturally  defined by the following formula

$$
G(\zeta,\eta):=\log{\vert\zeta\wedge\eta\vert\over\vert\zeta\vert\vert\eta\vert}=\log\sin{d(\zeta,\eta)\over\sqrt{2}}\,
\, \hbox{where} \,\, \vert\zeta\wedge\eta\vert^{2}=\sum_{0\leq i<j\leq n}\vert\zeta_{i}\eta_{j}-\zeta_{j}\eta_{i}\vert^{2} 
$$

and $d$ is the geodesic distance associated to the Fubini-Study metric (see \cite{15},\cite{6}). 

We recall some definitions and give a useful characterization of the local domain of definition of the complex Monge-Amp\`ere operator given by Z.
B\l ocki (see \cite{2},\cite{3}).

\begin{defi} Let $\Omega  \subset \C^n$  be a domain.   By definition the set $DMA_{loc}(\Omega)$ denotes the set of plurisubharmonic functions $\phi$ on $\Omega$ for which there a positive Borel measure $\sigma$ on $\Omega$ such that for all open $U\subset\subset X$ and $\forall\ (\phi_{j})\in PSH(U)\cap C^{\infty}(U)\searrow\phi$ in $U$, the sequences of measures $(dd^{c}\phi_{j})^{n}$ converges weakly to $\sigma$ in $U$. In this case, we put  $(dd^{c}\phi)^{n}=\sigma$.
\end{defi}

The following result of Blocki gives a useful characterization of the local domain of definition of the complex Monge-Amp\`ere operator. 
\begin{thm} (Z. B\l ocki \cite{2}, \cite{3}). 
1. If $\Omega\subset\mathbb{C}^{2}$ is an open set then $DMA_{loc}(\Omega)=PSH(\Omega)\cap W^{1,2}_{loc}(\Omega)$. 

2. If $n\geq 3$,  a plurisubharmonic function $\phi$ on a open set $U\subset\mathbb{C}^{n}$ belong to $DMA_{loc}(\Omega)$ if and only if  for any $z\in \Omega$ there exists a neighborhood $U_{z}$ of $z$  in $\Omega$  and a sequence $(\phi_{j})\subset PSH(U_{z})\cap C^{\infty}(U_{z})\searrow \phi$ in $U_{z}$ such that the sequences
$$
\vert \phi_{j}\vert^{n-p-2}d\phi_{j}\wedge d^{c}\phi_{j}\wedge(dd^{c}\phi_{j})^{p}\wedge(dd^{c}\vert z\vert^{2})^{n-p-1},\quad p=0,1,\cdots,n-2
$$
are locally weakly bounded in $U_{z}$.
\end{thm}

Observe that by Bedford and Taylor \cite{1}, the class of locally bounded plurisubharmonic functions in $\Omega$ is contained in $DMA_{loc}(\Omega)$.
By the work of J.-P. Demailly \cite{9}, any plurisubharmonic function in $\Omega$ bounded near the boundary $\partial \Omega$ is contained in $DMA_{loc}(\Omega)$.  
Let $(X,\omega)$ be a K\"ahler manifold of dimension $n$. We denote by $PSH (X,\omega)$  the set of $\omega$-plurisubharmonic functions in $X$. Then it is possible to define in the same way the local domain of definition
$DMA_{loc} (X,\omega)$ of the complex Monge-Amp\`ere operator on $(X,\omega)$.
A function $\varphi \in PSH (X,\omega)$ belongs to $DMA_{loc} (X,\omega)$  iff for any local chart $(U,z)$ the function $\phi := \varphi + \rho \in DMA_{loc} (U)$ where $\rho$ is a K\"ahler potential of $\omega$.
Then the previous theorem extends trivially to this general case.

Let $(\chi_{j})_{0\leq j\leq n}$ be a fixed partition of unity subordinated to the covering $(\mathcal{U}_{j})_{0\leq j\leq n}$. We define $m_{j}=\int\chi_{j}d\mu$ and $J=\{j\in\{0,1,\cdots,n\}\ :\ m_{j}\not=0\}$.  The $J\not=\emptyset$ and for $j\in J$, the measure $\mu_{j}:={1\over m_{j}}\chi_{j}\mu$ is a probability measure on $\mathbb{P}^{n}$ supported in $\mathcal{U}_{j}$ and we have the following convex decomposition of $\mu$
 $$
\mu=\sum_{j\in J}m_{j}\mu_{j}
 $$
 Therefore the potential $\mathbb{G}_{\mu}$ can be written as a convex combination
 $$
 \mathbb{G}_{\mu}=\sum_{j\in J} m_j \mathbb{G}_{\mu_{j}}.
 $$ 
 
 To show that $\mathbb{G}_{\mu}\in DMA_{loc}(\mathbb{P}^{n},\omega)$, it suffices to consider the case of a compact measure supported in an affine chart. Without loss of generality, we may always assume that $\mu$ is compactly supported in $\mathcal{U}_{0}$ and we are reduced to the study of the potential $\mathbb{G}_{\mu}$ on the open set $\mathcal{U}_{0}$. The restriction of $G(\zeta,\eta)$ to $\mathcal{U}_{0}\times\mathcal{U}_{0}$ can expressed in the affine coordinates as 
  $$
G(\zeta,\eta)=N(z,w)-{1\over 2}\log(1+\vert z\vert^{2})$$ where
$$
N(z,w):={1\over 2}\log{\vert z-w\vert^{2}+\vert z\wedge w\vert^{2}\over 1+\vert w\vert^{2}}
$$
 will be called the projective logarithmic kernel on $\mathbb{C}^{n}$.

\begin{lem}\label{diag}
1.The kernel $N$ is upper semicontinuous in $\mathbb{C}^{n}\times\mathbb{C}^{n}$ and smooth off the diagonal of $\mathbb{C}^{n}\times\mathbb{C}^{n}$.\\
2. For any fixed $w\in\mathbb{C}^{n}$, the function $ N(.,w)\ :\ z\rightarrow N(z,w)$ is plurisubharmonic in $\mathbb{C}^{n}$ and satisfies the following inequality
$$
{1\over 2}\log{\vert z-w\vert^{2}\over 1+\vert w\vert^{2}} \leq N(z,w)\leq{1\over 2}\log(1+\vert z\vert^{2}),\quad\forall\ (z,w)\in \mathbb{C}^{n}\times\mathbb{C}^{n}
$$
\end{lem}

From  lemma \ref{diag},  we have the following properties of the projective logarithmic kernel $G$ on $\mathbb{P}^{n}\times\mathbb{P}^{n}$.

\begin{cor}
1. The kernel $G$ is a non positive upper semi continuous function on $\mathbb{P}^{n}\times\mathbb{P}^{n}$ and smooth off the diagonal of $\mathbb{P}^{n}\times\mathbb{P}^{n}$.\\
2. For any fixed $\eta\in\mathbb{P}^{n}$, the function $G(.,\eta)\ :\ \zeta \rightarrow G(\zeta,\eta)$ is a non positive $\omega$-plurisubharmonic function in $\mathbb{P}^{n}$ and smooth in $\mathbb{P}^{n}\setminus\{\eta\}$, hence $G(.,\eta)\in  DMA_{loc}(\mathbb{P}^{n},\omega)$.  Moreover
$(\omega + dd^c G(\cdot,\eta))^n = \delta_\eta$.
\end{cor}

For a probability measure $\nu$ on $\mathbb{C}^{n}$, we define the projective logarithmic potential of $\nu$ as follows
$$
\mathbb{V}_{\nu}(z):={1\over 2}\int_{\mathbb{C}^{n}}\log{\vert z-w\vert^{2}+\vert z\wedge w\vert^{2}\over 1+\vert w\vert^{2}}d\nu(w)
$$

\begin{pro}\label{pluris}
 Let $\nu$ be a probability measure $\nu$ on $\mathbb{C}^{n}$. Then the function  $\mathbb{V}_{\nu}(z)$ is plurisubharmonic in $\mathbb{C}^{n}$ and  for all $z\in\mathbb{C}^{n}$
 $$
 \mathbb{V}_{\nu}(z)\leq{1\over 2}\log(1+\vert z\vert^{2}).
 $$

Also $\mathbb{V}_{\nu}\in DMA_{loc}(\mathbb{C}^{n})$ and
$$
(dd^{c}\mathbb{V}_{\nu})^{n}=\int_{\mathbb{C}^{n}\times\cdots\times\mathbb{C}^{n}}dd_{z}^{c}N(.,w_{1})\wedge\cdots\wedge dd_{z}^{c}N(.,w_{w})d\nu(w_{1})\cdots d\nu(w_{n})
$$
\end{pro}

\noindent{\bf Proof of theorem \ref{thm1}} As we have seen we have
$$
\G_\mu = \sum_{j \in J}  m_j  \G_{\mu_j},
$$
where  $\mu_j$ is compactly supported in the affine chart $\mathcal U_j$. 

Observe that for a fixed $k$ one can write  on $\mathcal U_k$
$$
\mathbb{G}_{\mu_k}(\zeta)+{1/2}\log(1+\vert z\vert^{2})=  \mathbb{V}_{\mu_k}(z), \, \, \text{where} \, \, z :=  z^{k} (\zeta) \in \C^n,
$$ 
 which is plurisubharmonic in $\C^n$.  Hence $\mathbb{G}_{\mu}$ is $\omega$-plurisubharmonic in $\mathbb{P}^{n}$.

2. By the co-area formula ( see  \cite{10} )
\begin{eqnarray*}
\int_{\mathbb{P}^{n}}\mathbb{G}_{\mu}(\zeta)dV(\zeta)&=&\int_{0}^{\pi/\sqrt{2}}\log\sin{r\over\sqrt{2}}A(r)dr\\
&=&-{c_{n}\over\sqrt{2}n^{2}}
\end{eqnarray*}
where $A(r):=c_{n}\sin^{2n-2}(r/\sqrt{2})\sin(\sqrt{2}r)$ is the area of the sphere about $\eta$ and radius $r$ on $\mathbb{P}^{n}$ and $c_{n}$ is a numerical constant  ( see \cite{14} page 168 or \cite{11} lemma 5.6] ).\\
Let $p\geq 1$. Since $\vert\nabla d\vert_{\omega}=1$, also by the co-area formula
\begin{eqnarray*}
\int_{\mathbb{P}^{n}}\vert\nabla\mathbb{G}_{\mu}(\zeta)\vert^{p}dV(\zeta)&\leq&
\int_{\mathbb{P}^{n}}\cot^{p}\Bigl({d(\zeta,\eta)\over\sqrt{2}}\Bigr)d\mu(\eta)dV(\zeta)\\&\leq& 2\sqrt{2}c_{n}\int_{0}^{\pi/2}\sin^{2n-1-p}t dt
\end{eqnarray*}
which is finite if and only if $p < 2n$. Hence for all $p\in ]0,2n[\ :\ \mathbb{G}_{\mu}\in W^{1,p}(\mathbb{P}^{n})$ ( by concavity of $x^{p}$).

3. When $n = 2$, we can apply the previous result to conclude that $ \mathbb{G}_{\mu}\in DMA_{loc}(\mathbb{P}^{2})$. When $n\geq 3$, we  apply Blocki's characterization staded above to show that $\mathbb{G}_{\mu_k} \in DMA_{loc}(\mathcal U_k)$.  We consider the  following approximating sequence
 $$
 \mathbb{V}^{\epsilon}_{\mu}(z)={1\over 2}\int_{\mathbb{C}^{n}}\log\Bigl({\vert z-w\vert^{2}+\vert z\wedge w\vert^{2}\over 1+\vert w\vert^{2}}+\epsilon^{2}\Bigr)d \mu_k(w)\searrow \mathbb{V}_{\mu}(z),
 $$
 and use the  next classical lemma on Riesz potentials to show a uniform estimates on their weighted  gradients as required in Blocki's theorem.
 
\begin{lem}
Let $\mu$ be a probability measure  on $\mathbb{C}^{n}$. For $0<\alpha<2n$, define the Riesz potential of $\mu$ by
$$
J_{\mu,\alpha}(z):=\int_{\mathbb{C}^{n}}{d\mu(w)\over\vert z-w\vert^{\alpha}}
$$
 If $0<p<2n/\alpha$ then $J_{\mu,\alpha}\in L^{p}_{loc}(\mathbb{C}^{n})$.
\end{lem}

\section{Regularizing property and proof of theorem \ref{thm2}}

We prove a regularizing property of the operator $\mu\rightarrow\mathbb{G}_{\mu}$. By the localization process exlained before, the proof of theorem 1.2 follows from the following theorem which generalizes and improves a result of Carlehed ( see \cite{5}).
\begin{thm}\label{thm3}
Let $\mu$ be a probability measure on $\mathbb{C}^{n}$  with no atoms and let $\psi \in \mathcal L (\C^n)$. Assume that $\psi$ is smooth in some open subset $U \subset \C^n$.  Then for any $0 \leq m \leq n $,  the Monge-Amp\`ere measure $(dd^{c}\mathbb{V}_{\mu})^{m} \wedge (dd^c \psi)^{n - m}$ is absolutely continuous with respect to the Lebesgue measure on $U$.
\end{thm}

The proof is based on the following elementary lemma.
\begin{lem}\label{absolut}
Let $(w_1, \cdots , w_n) \in (\C^n)^n$ fixed such that $w_{1}\not=w_{2}$. Let $\psi \in \mathcal L (\C^n)$. Assume that $\psi$ is smooth  in some open subset $U \subset \C^n$. Then for  any integer $ 0 \leq m \leq n$,   the measure 
$$
\bigwedge_{1 \leq j \leq m} dd^{c}\log(\vert \cdot -w_j\vert^{2}+\vert \cdot \wedge w_j\vert^{2}) \wedge(dd^c \psi) ^{n - m}
$$
is absolutely continuous with respect to the Lebesgue measure on $U$.

\end{lem}

{\bf Proof of theorem \ref{thm3}:} 
We first assume  that $m = n$. Let $K\subset\mathbb{C}^{n}$ be a compact set such that $(dd^{c}\vert z\vert^{2})^{n}(K)=0$. 
 Set $\Delta=\{(w,w,\cdots,w)\ :\ w\in\mathbb{C}^{n}\}$. Since $\mu$ puts no mass at any point, it follows by Fubini's theorem that $\mu^{\otimes n} (\Delta)=0$. By proposition \ref{pluris}
$$
\int_K (dd^{c}\mathbb{V}_{\mu})^n =\int_{(\mathbb{C}^{n})^{n}\setminus\Delta}f(w_{1},\cdots,w_{n})d\mu^{\otimes n}(w_{1},\cdots,w_{n})
$$
where 
$$
f(w_{1},\cdots,w_{n})=\int_{K}dd^{c}\log(\vert z-w_1\vert^{2}+\vert z\wedge w_1\vert^{2})\wedge\cdots\wedge d^{c}\log(\vert z-w_{n}\vert^{2}+\vert z\wedge w_{n}\vert^{2})
$$
By  Lemma \ref{absolut}, for any $(w_{1},\cdots,w_{n})\not\in\Delta$, $f(w_1, \cdots, w_n) = 0$, hence $(dd^{c}\mathbb{V}_{\mu})^{n}(K)=0$. The case $1 \leq m \leq n$ follows from Lemma \ref{absolut} in the same way.  The proof is complete.

\bigskip

{\bf Proof of theorem  \ref{thm2}:}  As we have seen in the proof of Theorem \ref{thm1}, one can write on each coordinate chart $\mathcal U_k$,
$$
\G_\mu (\zeta) = m_k \mathbb{G}_{\mu_k} + \psi_k (z),
$$
where $\psi_k \in \mathcal L (\C^n)$ is a smooth function in $\C^n$. Using Theorem \ref{thm3} again we conclude that $\G_\mu \in DMA_{loc} (\mathcal U_k)$.
Therefore $\G_\mu \in DMA_{loc} (\mathbb P^n)$.

% \section{}
% \label{}

% The Acknowledgements are an un-numbered section
%\section*{Acknowledgements}
% Acknowledgements text here

\section*{Acknowledgements}
It is a pleasure  to thank my supervisors  Ahmed Zeriahi and Said Asserda  for their support, suggestions and encouragement. I also would like to  thank professor Vincent Guedj for very useful discussions and suggestions. A part of this work was done when the author was visiting l'Institut de Math\'ematiques de Toulouse in March 2016. She would like to thank this institution for the invitation.

\end{document}